\input amstex
\documentstyle{amsppt}
%
\catcode`@=11
\redefine\output@{%
  \def\break{\penalty-\@M}\let\par\endgraf
  \ifodd\pageno\global\hoffset=105pt\else\global\hoffset=8pt\fi  
  \shipout\vbox{%
    \ifplain@
      \let\makeheadline\relax \let\makefootline\relax
    \else
      \iffirstpage@ \global\firstpage@false
        \let\rightheadline\frheadline
        \let\leftheadline\flheadline
      \else
        \ifrunheads@ 
        \else \let\makeheadline\relax
        \fi
      \fi
    \fi
    \makeheadline \pagebody \makefootline}%
  \advancepageno \ifnum\outputpenalty>-\@MM\else\dosupereject\fi
}
\def\Beta{\mathchar"0\hexnumber@\rmfam 42}
\catcode`\@=\active
\nopagenumbers
\chardef\textvolna='176
\def\negskp{\hskip -2pt}

\chardef\degree="5E
\def\blue#1{#1}

\catcode`#=11\def\diez{#}\catcode`#=6
\catcode`&=11\catcode`&=4
\catcode`_=11\def\podcherkivanie{_}\catcode`_=8
\catcode`~=11\def\volna{~}\catcode`~=\active
\def\mycite#1{\cite{\blue{#1}}\immediate\special{ps:
     ShrHPSdict begin /ShrBORDERthickness 0 def}}
\def\myciterange#1#2#3#4{\cite{\blue{#2#3#4}}\immediate\special{ps:
     ShrHPSdict begin /ShrBORDERthickness 0 def}}
\def\mytag#1{%
    \tag#1}
\def\mythetag#1{\thetag{\blue{#1}}\immediate\special{ps:
     ShrHPSdict begin /ShrBORDERthickness 0 def}}
\def\myrefno#1{\no#1}
\def\myhref#1#2{\blue{#2}\immediate\special{ps:
     ShrHPSdict begin /ShrBORDERthickness 0 def}}
\def\myEarXivlink{\myhref{http://arXiv.org}{http:/\negskp/arXiv.org}}

\def\mytheorem#1{\csname proclaim\endcsname{Theorem #1}}
\def\mytheoremwithtitle#1#2{\csname proclaim\endcsname{Theorem #1#2}}
\def\mythetheorem#1{\blue{#1}\immediate\special{ps:
     ShrHPSdict begin /ShrBORDERthickness 0 def}}
\def\mylemma#1{\csname proclaim\endcsname{Lemma #1}}
\def\mylemmawithtitle#1#2{\csname proclaim\endcsname{Lemma #1#2}}
\def\mythelemma#1{\blue{#1}\immediate\special{ps:
     ShrHPSdict begin /ShrBORDERthickness 0 def}}
\def\mycorollary#1{\csname proclaim\endcsname{Corollary #1}}

\def\myconjecture#1{\csname proclaim\endcsname{Conjecture #1}}
\def\myconjecturewithtitle#1#2{\csname proclaim\endcsname{Conjecture #1#2}}

\def\myproblem#1{\csname proclaim\endcsname{Problem #1}}
\def\myproblemwithtitle#1#2{\csname proclaim\endcsname{Problem #1#2}}
\def\mytheproblem#1{\blue{#1}\immediate\special{ps:
     ShrHPSdict begin /ShrBORDERthickness 0 def}}

\font\eightcyr=wncyr8
\pagewidth{360pt}
\pageheight{606pt}
\topmatter
\title
A note on solutions\\
of the cuboid factor equations.
\endtitle
\author
Ruslan Sharipov
\endauthor
\address Bashkir State University, 32 Zaki Validi street, 450074 Ufa, Russia
\endaddress
\email\myhref{mailto:r-sharipov\@mail.ru}{r-sharipov\@mail.ru}
\endemail
\abstract
     A rational perfect cuboid is a rectangular parallelepiped whose edges and 
face diagonals are given by rational numbers and whose space diagonal is equal 
to unity. It is described by a system of four quadratic equations with respect 
to six variables. The cuboid factor equations were derived from these four
equations by symmetrization procedure. They constitute a system of eight 
polynomial equations. Recently two sets of formulas were derived providing two
solutions for the cuboid factor equations. These two solutions are studied in 
the present paper. They are proved to coincide with each other up to a change 
of parameters in them.
\endabstract
\endtopmatter
\TagsOnRight
\document

\head
1. Introduction.
\endhead
     Finding a rational perfect cuboid is equivalent to finding a perfect
cuboid with all integer edges and diagonals, which is an old unsolved problem
known since 1719. The history of cuboid studies can be followed through the 
references \myciterange{1}{1}{--}{44}. Here are the equations describing 
perfect cuboids:
$$
\xalignat 2
&\hskip -2em
x_1^2+x_2^2+x_3^2-L^2=0,
&&x_2^2+x_3^2-d_1^{\kern 1pt 2}=0,\\
\vspace{-1.7ex}
\mytag{1.1}\\
\vspace{-1.7ex}
&\hskip -2em
x_3^2+x_1^2-d_2^{\kern 1pt 2}=0,
&&x_1^2+x_2^2-d_3^{\kern 1pt 2}=0.
\endxalignat
$$
The variables $x_1$, $x_2$, $x_3$ in \mythetag{1.1} represent edges of a cuboid,
the variables $d_1$, $d_2$, $d_3$ are its face diagonals, and $L$ is its space
diagonal. In the case of a rational perfect cuboid we set $L=1$.\par
     Let's denote through $p_{\kern 1pt 0}$, $p_{\kern 1pt 1}$, $p_{\kern 1pt 2}$, 
$p_{\kern 1pt 3}$ the left hand sides of the cuboid equations \mythetag{1.1}, i\.\,
e\. let's introduce the following notations:
$$
\xalignat 2
&\hskip -2em
p_{\kern 1pt 0}=x_1^2+x_2^2+x_3^2-L^2,
&&p_{\kern 1pt 1}=x_2^2+x_3^2-d_1^{\kern 1pt 2},\\
\vspace{-1.7ex}
\mytag{1.2}\\
\vspace{-1.7ex}
&\hskip -2em
p_{\kern 1pt 2}=x_3^2+x_1^2-d_2^{\kern 1pt 2},
&&p_{\kern 1pt 3}=x_1^2+x_2^2-d_3^{\kern 1pt 2}.
\endxalignat
$$
Using the polynomials \mythetag{1.2}, the following eight equations are written:
$$
\xalignat 2
&\hskip -2em
p_{\kern 1pt 0}=0,
&&\sum^3_{i=1}p_{\kern 1pt i}=0,\\
\vspace{-1.7ex}
\mytag{1.3}\\
\vspace{-1.7ex}
&\hskip -2em
\sum^3_{i=1}d_i\,p_{\kern 1pt i}=0,
&&\sum^3_{i=1}x_i\,p_{\kern 1pt i}=0,\\
\displaybreak
&\hskip -2em
\sum^3_{i=1}d_i^{\kern 1pt 2}\,p_{\kern 1pt i}=0,
&&\sum^3_{i=1}x_i^2\,p_{\kern 1pt i}=0,\\
\vspace{-1.7ex}
\mytag{1.4}\\
\vspace{-1.7ex}
&\hskip -2em
\sum^3_{i=1}x_i\,d_i\,p_{\kern 1pt i}=0,
&&\sum^3_{i=1}x_i^2\,d_i^{\kern 1pt 2}\,p_{\kern 1pt i}=0.\\
\endxalignat
$$
The equations \mythetag{1.3} and \mythetag{1.4} are called the cuboid 
factor equations. They were derived as a result of a symmetry approach 
to the original cuboid equations \mythetag{1.1} initiated in \mycite{45} 
(see also \myciterange{46}{46}{--}{48}).\par
     It is easy to see that each solution of the original cuboid equations 
\mythetag{1.1} is a solution for the factor equations \mythetag{1.3} and 
\mythetag{1.4}. Generally speaking, the converse is not true. However, in
\mycite{47} the following theorem was proved. 
\mytheorem{1.1} Each integer or rational solution of the factor equations 
\mythetag{1.3} and \mythetag{1.4} such that $x_1>0$, $x_2>0$, $x_3>0$, 
$d_1>0$, $d_2>0$, and $d_3>0$ is an integer or rational solution for the 
equations \mythetag{1.1}.
\endproclaim
      Due to the theorem~\mythetheorem{1.1} the factor equations \mythetag{1.3} 
and \mythetag{1.4} are equivalent to the equations \mythetag{1.1} in studying
perfect cuboids. But in this paper, saying a solution of the factor equations
we assume any integer or rational solution, i\.\,e\. even such that some of the 
inequalities $x_1>0$, $x_2>0$, $x_3>0$, $d_1>0$, $d_2>0$, $d_3>0$ or all
of them are not fulfilled.\par 
      Note that the left hand sides of the factor equations are multisimmetric
polynomials in $x_1$, $x_2$, $x_3$ and $d_1$, $d_2$, $d_3$, i\.\,e\. they are 
invariant with respect to the $S_3$ permutation group acting upon 
$x_1$, $x_2$, $x_3$, $d_1$, $d_2$, $d_3$, and $L$ as follows:
$$
\xalignat 3
&\sigma(x_i)=x_{\sigma i},
&&\sigma(d_i)=d_{\sigma i},
&&\sigma(L)=L.
\endxalignat
$$ 
For the theory of multisymmetric polynomials the reader is referred to
\myciterange{49}{49}{--}{69}. According to this theory, each multisymmetric
polynomial is expressed through the following nine elementary multisymmetric 
polynomials:
$$
\gather
\hskip -2em
\aligned
&x_1+x_2+x_3=E_{10},\\
&x_1\,x_2+x_2\,x_3+x_3\,x_1=E_{20},\\
&x_1\,x_2\,x_3=E_{30},
\endaligned
\mytag{1.5}\\
\vspace{2ex}
\hskip -2em
\aligned
&d_1+d_2+d_3=E_{01},\\
&d_1\,d_2+d_2\,d_3+d_3\,d_1=E_{02},\\
&d_1\,d_2\,d_3=E_{03},
\endaligned
\quad
\mytag{1.6}\\
\vspace{2ex}
\hskip -2em
\aligned
&x_1\,x_2\,d_3+x_2\,x_3\,d_1+x_3\,x_1\,d_2=E_{21},\\
&x_1\,d_2+d_1\,x_2+x_2\,d_3+d_2\,x_3+x_3\,d_1+d_3\,x_1=E_{11},\\
&x_1\,d_2\,d_3+x_2\,d_3\,d_1+x_3\,d_1\,d_2=E_{12}.
\endaligned
\mytag{1.7}
\endgather
$$
Expressing the left hand sides of the factor equations \mythetag{1.3} 
and \mythetag{1.4} through the polynomials \mythetag{1.5}, \mythetag{1.6}
and \mythetag{1.7}, one gets polynomial equations with respect to the 
variables $E_{10}$, $E_{20}$, $E_{30}$, $E_{01}$, $E_{02}$, $E_{03}$, 
$E_{21}$, $E_{11}$, $E_{12}$, and $L$ (see \thetag{3.1} through 
\thetag{3.7} in \mycite{48}). These equations were complemented with 
fourteen identities expressing the algebraic dependence of the 
elementary multisymmetric polynomials \mythetag{1.5}, \mythetag{1.6},
and \mythetag{1.7} (see \thetag{3.8} in \mycite{48}). As a result a 
system of twenty two polynomial equations was obtained. In \mycite{48} 
this huge system of twenty two polynomial equations was reduced 
to the following single polynomial equation for $E_{10}$, $E_{01}$,
$E_{11}$, and $L$:
$$
\hskip -2em
(2\,E_{11})^2+(E_{01}^2+L^2-E_{10}^2)^2-8\,E_{01}^2\,L^2=0. 
\mytag{1.8}
$$
The other variables $E_{20}$, $E_{30}$, $E_{02}$, $E_{03}$, $E_{21}$, 
$E_{12}$ are expressed as rational functions of $E_{10}$, $E_{01}$,
$E_{11}$, and $L$ (see formulas \thetag{4.1}, \thetag{4.3}, \thetag{5.1}, 
\thetag{5.2}, \thetag{4.6}, \thetag{4.7} in \mycite{48}).\par
     The equation \mythetag{1.8} was solved by John Ramsden in \mycite{70}. 
In the case of a rational perfect cuboid, where $L=1$, omitting some inessential
special cases, the general solution of the equation \mythetag{1.8} is given
by the formulas 
$$
\gather
\hskip -2em
E_{11}=-\frac{b\,(c^2+2-4\,c)}{b^2\,c^2+2\,b^2-3\,b^2\,c+c-b\,c^2\,+2\,b},
\mytag{1.9}\\
\vspace{1ex}
\hskip -2em
E_{10}=-\frac{b^2\,c^2+2\,b^2-3\,b^2\,c\,-c}{b^2\,c^2+2\,b^2-3\,b^2\,c
+c-b\,c^2+2\,b},
	\mytag{1.10}\\
\vspace{1ex}
\hskip -2em
E_{01}=-\frac{b\,(c^2+2-2\,c)}{b^2\,c^2+2\,b^2-3\,b^2\,c+c-b\,c^2+2\,b}.
\mytag{1.11}
\endgather
$$
Below are the formulas for $E_{12}$, $E_{21}$, $E_{03}$, $E_{30}$, $E_{02}$, 
$E_{20}$ in \mythetag{1.5}, \mythetag{1.6}, and \mythetag{1.7}:
$$
\allowdisplaybreaks
\gather
\hskip -2em
\gathered
E_{12}=(16\,b^6+32\,b^5-6\,c^5\,b^2+2\,c^5\,b-62\,b^5\,c^6
+62\,b^6\,c^6+16\,b^4\,-\\
-\,180\,b^6\,c^5-c^7\,b^3+18\,b^5\,c^7-12\,b^6\,c^7-2\,b^5\,c^8
+b^6\,c^8+248\,b^5\,c^2\,+\\
+\,248\,b^6\,c^2-96\,b^6\,c+321\,b^6\,c^4-180\,b^5\,c^3-144\,b^5\,c
-360\,b^6\,c^3\,+\\
+\,b^4\,c^8+8\,b^4\,c^6-6\,b^4\,c^7+18\,b^4\,c^5+7\,b^3\,c^6
+90\,b^5\,c^5-14\,b^3\,c^5\,+\\
+\,17\,b^2\,c^4+32\,b^4\,c^2+28\,b^3\,c^3-28\,b^3\,c^2-4\,b\,c^3+8\,b^3\,c
-57\,b^4\,c^4\,+\\
+\,36\,b^4\,c^3-12\,b^2\,c^3-48\,b^4\,c-c^4)\,(b^2\,c^4-6\,b^2\,c^3
+13\,b^2\,c^2\,-\\
-\,12\,b^2\,c+4\,b^2+c^2)^{-1}\,(b\,c-1-b)^{-2}\,(b\,c-c-2\,b)^{-2},
\endgathered\qquad
\mytag{1.12}\\
\vspace{2ex}
\gathered
E_{21}=\frac{b}{2}\,(5\,c^6\,b-2\,c^6\,b^2+52\,c^5\,b^2-16\,c^5\,b
-2\,c^7\,b^2+2\,b^4\,c^8\,-\\
-\,26\,b^4\,c^7-426\,b^4\,c^5-61\,b^3\,c^6+100\,b^3\,c^5
+14\,c^7\,b^3-c^8\,b^3-20\,b\,c^2\,-\\
-\,8\,b^2\,c^2-16\,b^2\,c-128\,b^2\,c^4-200\,b^3\,c^3
+244\,b^3\,c^2+32\,b\,c^3\,+\\
+\,768\,b^4\,c^4-852\,b^4\,c^3+568\,b^4\,c^2+104\,b^2\,c^3-208\,b^4\,c
+8\,c^4\,+\\
+16\,b^3-112\,b^3\,c+142\,b^4\,c^6
+32\,b^4-2\,c^5)\,(b^2\,c^4-6\,b^2\,c^3+13\,b^2\,c^2\,-\\
-12\,b^2\,c-4\,c^3+4\,b^2+c^2)^{-1}\,(b\,c-1-b)^{-2}\,(b\,c-c-2\,b)^{-2},
\endgathered\qquad\quad
\mytag{1.13}\\
\vspace{1ex}
\hskip -2em
\gathered
E_{03}=\frac{b}{2}\,(b^2\,c^4-5\,b^2\,c^3+10\,b^2\,c^2-10\,b^2\,c+4\,b^2
+2\,b\,c+2\,c^2\,-\\
-\,b\,c^3)\,(2\,b^2\,c^4-12\,b^2\,c^3+26\,b^2\,c^2-24\,b^2\,c
+\,8\,b^2-c^4\,b+3\,b\,c^3\,-\\
-\,6\,b\,c+4\,b+c^3-2\,c^2+2\,c)\,(b^2\,c^4-6\,b^2\,c^3+13\,b^2\,c^2\,-\\
-12\,b^2\,c+4\,b^2+c^2)^{-1}\,(b\,c-1-b)^{-2}\,(-c+b\,c-2\,b)^{-2},
\endgathered\qquad\quad
\mytag{1.14}\\
\vspace{1ex}
\hskip -2em
\gathered
E_{30}=c\,b^2\,(1-c)\,(c-2)\,(b\,c^2-4\,b\,c+2+4\,b)
\,(2\,b\,c^2-c^2-4\,b\,c\,+\\
+\,2\,b)\,(b^2\,c^4-6\,b^2\,c^3+13\,b^2\,c^2-12\,b^2\,c+4\,b^2
+c^2)^{-1}\,\times\\
\times\,(b\,c-1-b)^{-2}\,(-c+b\,c-2\,b)^{-2},
\endgathered\qquad\quad
\mytag{1.15}\\
\vspace{1ex}
\hskip -2em
\gathered
E_{02}=\frac{1}{2}\,(28\,b^2\,c^2-16\,b^2\,c-2\,c^2-4\,b^2-b^2\,c^4
+4\,b^3\,c^4-12\,b^3\,c^3\,+\\
+\,4\,b\,c^3+24\,b^3\,c-8\,b\,c-2\,b^4\,c^4+12\,b^4\,c^3-26\,b^4\,c^2
-8\,b^2\,c^3\,+\\
+24\,b^4\,c-16\,b^3-8\,b^4)\,(b\,c-1-b)^{-2}\,(b\,c-c-2\,b)^{-2},
\endgathered\qquad\quad
\mytag{1.16}\\
\vspace{1ex}
\hskip -2em
\gathered
E_{20}=\frac{b}{2}\,(b\,c^2-2\,c-2\,b)\,(2\,b\,c^2-c^2-6\,b\,c+2
+4\,b)\,\times\\
\times\,(b\,c-1-b)^{-2}\,(b\,c-c-2\,b)^{-2}.
\endgathered\qquad\quad
\mytag{1.17}
\endgather
$$
The formulas \mythetag{1.12}, \mythetag{1.13}, \mythetag{1.14}, 
\mythetag{1.15}, \mythetag{1.16}, \mythetag{1.17} were derived in 
\mycite{71} by substituting \mythetag{1.9}, \mythetag{1.10}, and 
\mythetag{1.11} along with $L=1$ into the corresponding formulas 
from \mycite{48}.\par
     Thus, the right hand sides of the equalities \mythetag{1.5},
\mythetag{1.6}, and \mythetag{1.7} turned out to be expressed through
two arbitrary rational parameters $b$ and $c$. The next step was to
resolve these equalities with respect to $x_1$, $x_2$, $x_3$, $d_1$, 
$d_2$, $d_3$. For this purpose in \mycite{71} the following two cubic 
equations were written:
$$
\align
&\hskip -2em
x^3-E_{10}\,x^2+E_{20}\,x-E_{30}=0,
\mytag{1.18}\\
\vspace{1ex}
&\hskip -2em
d^{\kern 1pt 3}-E_{01}\,d^{\kern 1pt 2}+E_{02}\,d-E_{03}=0.
\mytag{1.19}
\endalign
$$
Note that the left hand sides of the equalities \mythetag{1.5} are regular
symmetric polynomials of the variables $x_1$, $x_2$, $x_3$ (see \mycite{72}). 
Similarly, the left hand sides of the equalities \mythetag{1.6} are regular
symmetric polynomials of the variables $d_1$, $d_2$, $d_3$. For this reason
$x_1$, $x_2$, $x_3$ can be found as roots of the cubic equation 
\mythetag{1.18}. Similarly, $d_1$, $d_2$, $d_3$ are roots of the second
cubic equation \mythetag{1.19}. Relying on these facts, in \mycite{71} the
following two inverse problems were formulated. 
\myproblem{1.1} Find all pairs of rational numbers $b$ and $c$ for which the
cubic equations \mythetag{1.18} and \mythetag{1.19} with the coefficients given
by the formulas \mythetag{1.10}, \mythetag{1.11},	\mythetag{1.14}, \mythetag{1.15}, 
\mythetag{1.16},	\mythetag{1.17} possess positive rational roots $x_1$, $x_2$, 
$x_3$, $d_1$, $d_2$, $d_3$ obeying the auxiliary polynomial equations 
\mythetag{1.7} whose right hand sides are given by the formulas \mythetag{1.9}, 
\mythetag{1.12}, and \mythetag{1.13}. 
\endproclaim
\myproblem{1.2} Find at least one pair of rational numbers $b$ and $c$ for which 
the cubic equations \mythetag{1.18} and \mythetag{1.19} with the coefficients given
by the formulas \mythetag{1.10}, \mythetag{1.11},	\mythetag{1.14}, \mythetag{1.15}, 
\mythetag{1.16},	\mythetag{1.17} possess positive rational roots $x_1$, $x_2$, 
$x_3$, $d_1$, $d_2$, $d_3$ obeying the auxiliary polynomial equations 
\mythetag{1.7} whose right hand sides are given by the formulas \mythetag{1.9}, 
\mythetag{1.12}, and \mythetag{1.13}.
\endproclaim
     Due to the theorem~\mythetheorem{1.1} the inverse problems~\mytheproblem{1.1} 
and \mytheproblem{1.2} are equivalent to finding all rational perfect cuboids and 
to finding at least one rational perfect cuboid respectively. Singularities of 
the inverse problems~\mytheproblem{1.1} and \mytheproblem{1.2} due to the 
denominators in the formulas \mythetag{1.9} through \mythetag{1.17} were studied
in \mycite{73}. Some special cases where the equations \mythetag{1.5},
\mythetag{1.6}, \mythetag{1.7} are solvable with respect to the cuboid variables
$x_1$, $x_2$, $x_3$ and $d_1$, $d_2$, $d_3$ were found in \mycite{74}. \pagebreak 
However, none of these special cases have produced a perfect cuboid since the 
inequalities 
$$
\xalignat 3
&\hskip -2em
x_1>0, &&x_2>0, &&x_3>0,\\
\vspace{-1.7ex}
\mytag{1.20}\\
\vspace{-1.7ex}
&\hskip -2em
d_1>0, &&d_2>0, &&d_3>0
\endxalignat 
$$
required for solving the problems~\mytheproblem{1.1} and \mytheproblem{1.2} are
not fulfilled in these special cases.\par
     Again, neglecting the inequalities \mythetag{1.20}, an approach to solving
the equations \mythetag{1.5}, \mythetag{1.6}, \mythetag{1.7} was found in 
\mycite{75}. It exploits the following lemma. 
\mylemma{1.1} A reduced cubic equation $y^3+y^2+D=0$ has three rational roots 
if and only if there is a rational number $w$ satisfying the sextic equation 
$$
\hskip -2em
D\,(w^2+3)^3+4\,(w-1)^2\,(1+w)^2=0.
\mytag{1.21}
$$
In this case the roots of the cubic equation $y^3+y^2+D=0$ are given by the
formulas 
$$
\xalignat 3
&\hskip -2em
y_1=-\frac{2\,(w+1)}{w^2+3},   
&&y_2=\frac{2\,(w-1)}{w^2+3},
&&y_3=\frac{1-w^2}{w^2+3}.
\quad
\mytag{1.22}
\endxalignat
$$
\endproclaim
     Based on the lemma~\mythelemma{1.1} and on the cubic equations \mythetag{1.18} 
and \mythetag{1.19}, in \mycite{75} two sextic equations of the form \mythetag{1.21} 
were derived:
$$
\align
&\hskip -2em
D_1\,(w_1^2+3)^3+4\,(w_1-1)^2\,(1+w_1)^2=0,
\mytag{1.23}\\
&\hskip -2em
D_2\,(w_2^2+3)^3+4\,(w_2-1)^2\,(1+w_2)^2=0.
\mytag{1.24}\\
\endalign
$$
The $D$-parameters $D_1$ and $D_2$ of the sextic equations \mythetag{1.23} and
\mythetag{1.24} depend on the same two rational numbers $b$ and $c$ as $E_{11}$,  
$E_{10}$, $E_{01}$, $E_{12}$, $E_{21}$, $E_{03}$, $E_{30}$, $E_{02}$, $E_{20}$ 
in the formulas \mythetag{1.9} through \mythetag{1.17}. They are given
by the formulas
$$
\gathered
D_1=-\frac{2}{27}\,(7812\,b^4\,c^4\,-216\,b^2\,c^4-52\,b^2\,c^3+1764\,b^3\,c^4
-1200\,b^4\,c^3\,-\\
-\,1848\,b^4\,c^2+720\,b^4\,c-36\,c^4\,b-1512\,b^3\,c^3-36\,c^8\,b^3
+288\,b^3\,c^2\,-\\
-\,108\,c^6\,b^2+380\,c^5\,b^2+378\,c^7\,b^3-231\,c^8\,b^4-300\,c^7\,b^4
+3906\,c^6\,b^4\,-\\
-13\,c^7\,b^2-8904\,c^5\,b^4-882\,c^6\,b^3+18\,c^6\,b-1319\,b^6\,c^8
+20952\,b^5\,c^3\,-\\
-\,11952\,b^5\,c^2+2592\,b^5\,c-48372\,b^6\,c^4+31620\,b^6\,c^3-10552\,b^6\,c^2\,+\\
+\,\,816\,b^6\,c+1494\,b^5\,c^8-5238\,b^5\,c^7-4\,c^5+7905\,b^6\,c^7
-24186\,b^6\,c^6\,+\\
+\,288\,b^6+43740\,b^6\,c^5+7686\,b^5\,c^6+576\,b^7+128\,b^8-15372\,b^5\,c^4\,-\\
-\,1080\,b^7\,c^8-3546\,b^7\,c^6+51\,c^9\,b^6+400\,b^8\,c^8-162\,c^9\,b^5
+8640\,b^7\,c^2\,-\\
-\,3456\,b^7\,c+2808\,b^7\,c^7-1560\,b^8\,c^7+3940\,b^8\,c^6+216\,c^9\,b^7
-960\,b^8\,c\,-\\
-\,6240\,b^8\,c^3+9\,c^{10}\,b^6+7880\,b^8\,c^4+4\,c^{10}\,b^8-6732\,b^8\,c^5 
+45\,c^9\,b^4\,+\\
+\,3200\,b^8\,c^2-11232\,b^7\,c^3+7092\,b^7\,c^4-18\,c^{10}\,b^7
-60\,c^9\,b^8)^2\,(2\,c^2\,+\\
+\,2\,b^4\,c^4-12\,b^4\,c^3+26\,b^4\,c^2-24\,b^4\,c+8\,b^4-6\,b^3\,c^4
+18\,b^3\,c^3\,-\\
-\,36\,b^3\,c+24\,b^3+3\,b^2\,c^4+8\,b^2\,c^3-36\,b^2\,c^2+16\,b^2\,c+12\,b^2
-6\,b\,c^3\,+\\
+\,12\,b\,c)^{-3}\,(b^2\,c^4-6\,b^2\,c^{-3}+13\,b^2\,c^2-12\,b^2\,c+4\,b^2+c^2)^{-2},
\endgathered
\mytag{1.25}
$$
$$
\gathered
D_2=-\frac{2\,b^2}{27}\,(832\,b^2\,c^2-1440\,b^2\,c^4-840\,b^2\,c^3
+4788\,b^3\,c^4+396\,b\,c^3\,+\\
+\,720\,b^3\,c+808\,b^4\,c^4+3032\,b^4\,c^3-2576\,b^4\,c^2
-96\,b^4\,c+448\,b^4\,-\\
-\,504\,c^4\,b-4176\,b^3\,c^3-9\,c^8\,b^3+72\,b^3\,c^2
-720\,c^6\,b^2+2288\,c^5\,b^2\,+\\
+\,1044\,c^7\,b^3-322\,c^8\,b^4+758\,c^7\,b^4+404\,c^6\,b^4
-210\,c^7\,b^2-2464\,c^5\,b^4\,-\\
-\,2394\,c^6\,b^3+72\,c^4+252\,c^6\,b+3168\,b^6\,c^8+441\,c^9\,b^5
-7056\,b^5\,c\,+\\
+\,57960\,b^6\,c^4-47232\,b^6\,c^3+25344\,b^6\,c^2
-8064\,b^6\,c-1809\,b^5\,c^8\,+\\
+\,14472\,b^5\,c^2+3951\,b^5\,c^7-72\,c^5+36\,c^6-11808\,b^6\,c^7
+1440\,b^5\,+\\
+\,28980\,b^6\,c^6-49032\,b^6\,c^5-4410\,b^5\,c^6
+8820\,b^5\,c^4-15804\,b^5\,c^3\,+\\
+\,1152\,b^6-504\,c^9\,b^6-45\,c^9\,b^3-6\,c^9\,b^4+104\,c^8\,b^2
+36\,c^{10}\,b^6\,+\\
+\,14\,c^{10}\,b^4-45\,c^{10}\,b^5-99\,c^7\,b)^2\,(6\,b^4\,c^4
-36\,b^4\,c^3+78\,b^4\,c^2
-72\,b^4\,c\,+\\
+\,24\,b^4-12\,b^3\,c^4+36\,b^3\,c^3-72\,b^3\,c+48\,b^3+5\,b^2\,c^4
+16\,b^2\,c^3\,-\\
-\,68\,b^2\,c^2+32\,b^2\,c+20\,b^2-12\,b\,c^3+24\,b\,c
+6\,c^2)^{-3}\,(b^2\,c^4-6\,b^2\,c^3\,+\\
+\,13\,b^2\,c^2-12\,b^2\,c+4\,b^2+c^2)^{-2}.
\endgathered\quad
\mytag{1.26}
$$
Along with \mythetag{1.25} and \mythetag{1.26}, in \mycite{75} twelve 
rational functions were derived:
$$
\xalignat 3
&\hskip -2em
x_1=x_1(b,c,w_1), &&x_2=x_2(b,c,w_1), &&x_3=x_3(b,c,w_1),\qquad\\
\vspace{-1.7ex}
\mytag{1.27}\\
\vspace{-1.7ex}
&\hskip -2em
d_1=d_1(b,c,w_1), &&d_2=d_2(b,c,w_1), &&d_3=d_3(b,c,w_1),\qquad\\
\vspace{2ex}
&\hskip -2em
x_1=x_1(b,c,w_2), &&x_2=x_2(b,c,w_2), &&x_3=x_3(b,c,w_2),\qquad\\
\vspace{-1.7ex}
\mytag{1.28}\\
\vspace{-1.7ex}
&\hskip -2em
d_1=d_1(b,c,w_2), &&d_2=d_2(b,c,w_2), &&d_3=d_3(b,c,w_2).\qquad
\endxalignat
$$
The explicit formulas for \mythetag{1.27} and \mythetag{1.28} are 
very huge. Therefore we provide them in the ancillary files 
{\bf Solutions\_\kern 1.5pt 1.txt} and {\bf Solutions\_\kern 1.5pt 2.txt} 
attached to this arXiv submission. The main result of \mycite{75} is 
formulated in the following theorems.
\mytheorem{1.2} If\/ $(b,c,w_1)$ is a triple of rational numbers
solving the equation \mythetag{1.23}, where $D_1$ is given by 
\mythetag{1.25}, and belonging to the domain of the rational functions
\mythetag{1.27}, then the values of these functions provide a rational 
solution for the equations \mythetag{1.5}, \mythetag{1.6}, \mythetag{1.7} 
and for the cuboid factor equations \mythetag{1.3}, \mythetag{1.4}.
\endproclaim
\mytheorem{1.3} If\/ $(b,c,w_2)$ is a triple of rational numbers
solving the equation \mythetag{1.24}, where $D_2$ is given by 
\mythetag{1.26}, and belonging to the domain of the rational functions
\mythetag{1.28}, then the values of these functions provide
a rational solution for the equations \mythetag{1.5}, \mythetag{1.6}, 
\mythetag{1.7} and for the cuboid factor equations \mythetag{1.3}, 
\mythetag{1.4}.
\endproclaim
    The main goal of the present paper is to prove that the sets of
solutions to the equations \mythetag{1.5}, \mythetag{1.6}, \mythetag{1.7} 
and to the cuboid factor equations \mythetag{1.3} and \mythetag{1.4} provided
by the theorems~\mythetheorem{1.2} and \mythetheorem{1.3} do essentially 
coincide.
\par 
\head
2. Some prerequisites.
\endhead
     Let's consider a general cubic equation with the coefficients $A_0$, $A_1$, 
$A_2$, $A_3$:
$$
\hskip -2em
A_3\,x^3+A_2\,x^2+A_1\,x+A_0=0.
\mytag{2.1}
$$
Under some certain restrictions for its coefficients, the cubic equation
\mythetag{2.1} can be transformed to its reduced form $y^3+y^2+D=0$, where $D$
is given by the formula
$$
\hskip -2em
D=-\frac{(9\,A_1\,A_2\,A_3-27\,A_0\,A_3^2-2\,A_2^3)^2}
{27\,(A_2^2-3\,A_1\,A_3)^3}.
\mytag{2.2}
$$
Applying the lemma~\mythelemma{1.1} to the reduced form of the equation 
\mythetag{2.1}, one gets the formulas \mythetag{1.22} for $y_1$, $y_2$, $y_3$. 
Then the backward transformation of $y_1$, $y_2$, $y_3$ to the roots of the 
equation \mythetag{2.1} yields the formulas
$$
\gather
\hskip -2em
\gathered
x_1=\frac{1}{18}\,
((2\,A_2^3-9\,A_1\,A_2\,A_3+27\,A_0\,A_3^2)\,w^2
+(18\,A_2\,A_1\,A_3-6\,A_2^3)\,w\,-\\
-\,9\,A_1\,A_2\,A_3+81\,A_0\,A_3^2
)\,A_3^{-1}\,(A_2^2-3\,A_1\,A_3)^{-1}\,(1+w)^{-1},
\endgathered\quad
\mytag{2.3}\\
\hskip -2em
\gathered
x_2=\frac{1}{18}\,
((2\,A_2^3-9\,A_1\,A_2\,A_3+27\,A_0\,A_3^2)\,w^2
-(18\,A_2\,A_1\,A_3-6\,A_2^3)\,w\,-\\
-\,9\,A_1\,A_2\,A_3+81\,A_0\,A_3^2
)\,A_3^{-1}\,(A_2^2-3\,A_1\,A_3)^{-1}\,(1-w)^{-1},
\endgathered\quad
\mytag{2.4}\\
\vspace{2ex}
\gathered
x_3=\frac{1}{9}\,((A_2^3-27\,A_0\,A_3^2)\,w^2+36\,A_1\,A_2\,A_3
-81\,A_0\,A_3^2-9\,A_2^3)\,\times\\
\times\,A_3^{-1}\,(A_2^2-3\,A_1\,A_3)^{-1}\,(1-w)^{-1}\,(1+w)^{-1}.
\endgathered\qquad
\mytag{2.5}
\endgather
$$
As a result one can formulate the following lemma for the equation
\mythetag{2.1}.
\mylemma{2.1} Assume that the numbers $A_0$, $A_1$, $A_2$, $A_3$ obey the 
inequalities
$$
\xalignat 3
&\hskip -2em
A_3\neq 0, 
&&\frac{A_1}{A_3}-\frac{A_2^2}{3\,A_3^2}\neq 0,
&&\frac{A_0}{A_3}-\frac{A_1\,A_2}{3\,A_3^2}+\frac{2\,A_2^3}{27\,A_3^3}\neq 0.
\qquad
\mytag{2.6}
\qquad
\endxalignat
$$
Then the general cubic polynomial \mythetag{2.1} with the rational coefficients 
$A_0$, $A_1$, $A_2$, $A_3$ has three rational roots if and only if there is a 
rational number $w$ satisfying the sextic equation \mythetag{2.1} where $D$ is
given by the formula \mythetag{2.2}. In this case the roots of the cubic equation 
\mythetag{2.1} are given by the formulas \mythetag{2.3}, \mythetag{2.4},
\mythetag{2.5}. 
\endproclaim 
     The detailed proofs of the lemmas~\mythelemma{1.1} and \mythelemma{2.1} can
be found in \mycite{75}.\par
     Now, assume that we have a cubic equation with three rational roots $x_1$,
$x_2$, $x_3$. Then this cubic equation can be written as follows:
$$
\hskip -2em
(x-x_1)(x-x_2)(x-x_3)=0
\mytag{2.7}
$$
Expanding the left hand side of the equation \mythetag{2.7}, we find
$$
\xalignat 2
&\hskip -2em
A_3=1,
&&A_1=x_1\,x_2+x_2\,x_3+x_3\,x_1,\\
\vspace{-1.7ex}
\mytag{2.8}\\
\vspace{-1.7ex}
&\hskip -2em
A_0=-x_1\,x_2\,x_3,
&&A_2=-(x_1+x_2+x_3).
\endxalignat
$$
The condition $A_3\neq 0$ from \mythetag{2.6} is fulfilled for the polynomial
\mythetag{2.7} since $A_3=1$ in \mythetag{2.8}. The second condition 
\mythetag{2.6} for the polynomial \mythetag{2.7} is written as 
$$
\hskip -2em
x_1^2+x_2^2+x_3^2-x_2\,x_3-x_1\,x_3-x_1\,x_2\neq 0.
\mytag{2.9}
$$
The third condition \mythetag{2.6} is the most interesting of all three. 
\pagebreak Applying \mythetag{2.8} to it, we find that for the polynomial 
\mythetag{2.7} this condition is written as follows:
$$
\hskip -2em
(2\,x_1-x_2-x_3)\,(2\,x_2-x_3-x_1)\,(2\,x_3-x_1-x_2)\neq 0.
\mytag{2.10}
$$
The condition \mythetag{2.10} can be written as three conditions
$$
\align
&\hskip -2em
u_1=2\,x_1-x_2-x_3\neq 0,\\
&\hskip -2em
u_2=2\,x_2-x_3-x_1\neq 0,
\mytag{2.11}\\
&\hskip -2em
u_3=2\,x_3-x_1-x_2\neq 0.
\endalign
$$
It is not obvious, but the condition \mythetag{2.9} can be written as 
$$
(2\,x_1-x_2-x_3)^2+(2\,x_2-x_3-x_1)^2+(2\,x_3-x_1-x_2)^2\neq 0.
$$
Therefore it is clear that the conditions \mythetag{2.11} imply both 
\mythetag{2.9} and \mythetag{2.10}.\par
     Now let's substitute \mythetag{2.8} into \mythetag{2.2}. Then we find 
that the $D$-parameter of the sextic equation \mythetag{1.21} corresponding 
to the cubic equation \mythetag{2.7} is written as
$$
\hskip -2em
\gathered
D=-\frac{8\,(u_1\,u_2\,u_3)^2}{(u_1^2+u_2^2+u_3^2)^3\strut}
\endgathered
\mytag{2.12}
$$
The denominator of \mythetag{2.12} is nonzero due to \mythetag{2.11}. 
Substituting \mythetag{2.12} into the sextic equation, we find that
it factors explicitly
$$
\hskip -2em
D\,(w-\tilde w_1)\,(w-\tilde w_2)\,(w-\tilde w_3)
\,(w-\tilde w_4)\,(w-\tilde w_5)\,(w-\tilde w_6)=0,
\mytag{2.13}
$$
where $D$ is given by \mythetag{2.12} and $\tilde w_1$, $\tilde w_2$, 
$\tilde w_3$ $\tilde w_4$ $\tilde w_5$, $\tilde w_6$ are given by the 
formulas
$$
\xalignat 2
&\hskip -2em
\tilde w_1=\frac{u_1-u_2}{u_3},
&&\tilde w_2=-\frac{u_1-u_2}{u_3},\\
\vspace{1ex}
&\hskip -2em
\tilde w_3=\frac{u_2-u_3}{u_1},
&&\tilde w_4=-\frac{u_2-u_3}{u_1},
\mytag{2.14}\\
\vspace{1ex}
&\hskip -2em
\tilde w_5=\frac{u_3-u_1}{u_2},
&&\tilde w_6=-\frac{u_3-u_1}{u_2}.
\endxalignat
$$
The numbers $u_1$, $u_2$, $u_3$ in \mythetag{2.12} and \mythetag{2.14} are 
determined by the formulas \mythetag{2.11}.\par
     The quantities $\tilde w_1$, $\tilde w_2$, $\tilde w_3$, $\tilde w_4$, 
$\tilde w_5$, $\tilde w_6$ in \mythetag{2.13} are roots of the sextic equation 
\mythetag{1.21}. Therefore, substituting \mythetag{2.8} into \mythetag{2.3}, 
\mythetag{2.4}, \mythetag{2.5} and substituting one of the quantities 
\mythetag{2.14} for $w$ into these formulas, we express $x_1$, $x_2$, $x_3$ 
through $x_1$, $x_2$, $x_3$, but up to some permutation of them. Here are the 
permutations associated with the quantities $\tilde w_1$, $\tilde w_2$, 
$\tilde w_3$, $\tilde w_4$, $\tilde w_5$, $\tilde w_6$ from \mythetag{2.14}:
$$
\xalignat 2
&\hskip -2em
\tilde w_1\!:\,(x_1,x_2,x_3)\mapsto (x_1,x_2,x_3),
&&\tilde w_2\!:\,(x_1,x_2,x_3)\mapsto (x_2,x_1,x_3),\\
&\hskip -2em
\tilde w_3\!:\,(x_1,x_2,x_3)\mapsto (x_2,x_3,x_1),
&&\tilde w_4\!:\,(x_1,x_2,x_3)\mapsto (x_3,x_2,x_1),
\qquad
\mytag{2.15}\\
&\hskip -2em
\tilde w_5\!:\,(x_1,x_2,x_3)\mapsto (x_3,x_1,x_2),
&&\tilde w_6\!:\,(x_1,x_2,x_3)\mapsto (x_1,x_2,x_3).
\endxalignat
$$
As we see in \mythetag{2.15}, the first quantity $\tilde w_1$ from 
\mythetag{2.14} plays the role of the identical permutation belonging to the 
permutation group $S_3$ and being its unit element.\par
\head
3. Conversion formulas. 
\endhead
     Let's choose the first formula \mythetag{2.14}. It expresses a root of 
the sextic equation \mythetag{1.21} through the roots of the associated cubic
equation \mythetag{2.1} given by the formulas \mythetag{2.3}, \mythetag{2.4}, 
\mythetag{2.5}. We write this formula as follows:
$$
\hskip -2em
w=\frac{3\,(x_1-x_2)}{2\,x_3-x_1-x_2}.
\mytag{3.1}
$$
The formula \mythetag{3.1} is inverse to the formulas \mythetag{2.3}, 
\mythetag{2.4}, and \mythetag{2.5}. Indeed, applying the formulas \mythetag{2.8}
to \mythetag{2.3}, \mythetag{2.4}, and \mythetag{2.5}, we obtain
$$
\gather
\hskip -2em
\gathered
x_1=\frac{1}{18}\bigl((x_1+x_2-2\,x_3)\,(x_2+x_3-2\,x_1)\,(x_3+x_1
-2\,x_2)\,w^2\,+\\
+\,6\,(x_3+x_2+x_1)\,(x_1^2+x_2^2+x_3^2-x_2\,x_3-x_1\,x_3-x_1\,x_2)\,w\,+\\
+\,9\,x_1\,x_2^2+9\,x_1\,x_3^2+9\,x_2\,x_3^2+9\,x_2\,x_1^2
+9\,x_3\,x_1^2+9\,x_3\,x_2^2\,-\\
-\,54\,x_1\,x_2\,x_3\bigr)\,(x_1^2+x_2^2+x_3^2-x_1\,x_2-x_2\,x_3
-x_3\,x_1)^{-1}\,(1+w)^{-1},
\endgathered
\mytag{3.2}\\
\vspace{2ex}
\hskip -2em
\gathered
x_2=\frac{1}{18}\bigl((x_1+x_2-2\,x_3)\,(x_2+x_3-2\,x_1)\,(x_3+x_1
-2\,x_2)\,w^2\,-\\
-\,6\,(x_3+x_2+x_1)\,(x_1^2+x_2^2+x_3^2-x_2\,x_3-x_1\,x_3-x_1\,x_2)\,w\,+\\
+\,9\,x_1\,x_2^2+9\,x_1\,x_3^2+9\,x_2\,x_3^2+9\,x_2\,x_1^2
+9\,x_3\,x_1^2+9\,x_3\,x_2^2\,-\\
-\,54\,x_1\,x_2\,x_3\bigr)\,(x_1^2+x_2^2+x_3^2-x_1\,x_2-x_2\,x_3
-x_3\,x_1)^{-1}\,(1-w)^{-1},
\endgathered
\mytag{3.3}\\
\vspace{2ex}
\hskip -2em
\gathered
x_3=\frac{1}{9}\bigl((21\,x_1\,x_2\,x_3+3\,x_1\,x_2^2+3\,x_1\,x_3^2
+3\,x_2\,x_3^2+3\,x_2\,x_1^2+3\,x_3\,x_1^2\,+\\
+\,3\,x_3\,x_2^2-x_1^3-x_2^3-x_3^3)\,w^2-9\,x_1\,x_2^2
-9\,x_1\,x_3^2-9\,x_2\,x_3^2-9\,x_2\,x_1^2\,-\\
-\,9\,x_3\,x_1^2-9\,x_3\,x_2^2+27\,x_1\,x_2\,x_3+9\,x_1^3+9\,x_2^3
+9\,x_3^3\bigr)\,\times\\
\times\,(x_1^2+x_2^2+x_3^2-x_1\,x_2-x_2\,x_3
-x_3\,x_1)^{-1}\,(1-w)^{-1}\,(1+w)^{-1}.
\endgathered
\quad
\mytag{3.4}
\endgather
$$     
Substituting \mythetag{3.2}, \mythetag{3.3}, \mythetag{3.4} into \mythetag{3.1}, 
we get the identity $w=w$. Conversely, substituting \mythetag{3.1} into 
\mythetag{3.2}, \mythetag{3.3}, \mythetag{3.4}, we get three identities
$x_1=x_1$, $x_2=x_2$, and $x_3=x_3$. This result proves the following theorem. 
\mytheorem{3.1} Let\/ $x_1$, $x_2$, $x_3$ be three roots of a general cubic 
equation \mythetag{2.1} such that the conditions \mythetag{2.6} are fulfilled. 
Then the formula \mythetag{3.1} yields a solution $w$ of the associated sextic 
equation \mythetag{1.21} whose $D$-parameter is given by the formula \mythetag{2.2}. 
In this case the roots $x_1$, $x_2$, $x_3$ are backward expressed through $w$ 
by means of the formulas \mythetag{2.3}, \mythetag{2.4}, and \mythetag{2.5}. 
\endproclaim
     Let's return to the formulas \mythetag{1.27},  \mythetag{1.28} and let's recall 
that the functions $x_1(b,c,w_1)$, $x_2(b,c,w_1)$, $x_3(b,c,w_1)$ from \mythetag{1.27}
were produced in \mycite{75} by applying the formulas \mythetag{2.3}, \mythetag{2.4}, 
\mythetag{2.5} to the cubic equation \mythetag{1.18}. Therefore, setting $w=w_1$ and
substituting these functions for $x_1$, $x_2$, $x_3$ into \mythetag{3.1}, we get the 
identity $w_1=w_1$. However, we have the other three functions $x_1(b,c,w_2)$, 
$x_2(b,c,w_2)$, $x_3(b,c,w_2)$ in \mythetag{1.28}. They represent the same three 
roots $x_1$, $x_2$, $x_3$ of the same cubic equation \mythetag{1.18}. Substituting 
them into \mythetag{3.1}, \pagebreak we get the same quantity $w_1$. But now $w_1$ 
turns out to be expressed through $w_2$, i\.\,e\. we get the formula 
$$
\hskip -2em
w_1=\frac{3\,x_1(b,c,w_2)-3\,x_2(b,c,w_2)}{2\,x_3(b,c,w_2)-x_1(b,c,w_2)
-x_2(b,c,w_2)},
\mytag{3.5}
$$
which is not an identity. The formula \mythetag{3.5} is the first conversion formula. 
It expresses $w_1$ through $b$, $c$, and $w_2$, i\.\,e\. \mythetag{3.5} yields a 
function $w_1=w_1(b,c,w_2)$. If we substitute this function into the argument
$w_1$ of the functions $x_1(b,c,w_1)$, $x_2(b,c,w_1)$, $x_3(b,c,w_1)$ from 
\mythetag{1.27}, then, according to the theorem~\mythetheorem{3.1}, we get back
three roots $x_1(b,c,w_2)$, $x_2(b,c,w_2)$, $x_3(b,c,w_2)$ used in \mythetag{3.5}.
This means that we have the following relationships based on the formula
\mythetag{3.5}:
$$
\align
&\hskip -2em
x_1(b,c,w_1(b,c,w_2))=x_1(b,c,w_2),\\
&\hskip -2em
x_2(b,c,w_1(b,c,w_2))=x_2(b,c,w_2),
\mytag{3.6}\\
&\hskip -2em
x_3(b,c,w_1(b,c,w_2))=x_3(b,c,w_2).
\endalign	
$$\par
     Note that \mythetag{1.19} is another cubic equation. It is associated
with the other sextic equation \mythetag{1.24} and it has its own formula like 
\mythetag{3.1}: 
$$
\hskip -2em
w=\frac{3\,(d_1-d_2)}{2\,d_3-d_1-d_2}.
\mytag{3.7}
$$
The functions $d_1(b,c,w_2)$, $d_2(b,c,w_2)$, $d_3(b,c,w_2)$ were produced 
in \mycite{75} by applying the formulas \mythetag{2.3}, \mythetag{2.4}, 
\mythetag{2.5} to the cubic equation \mythetag{1.19}. Therefore, setting $w=w_2$ 
and substituting these functions for $d_1$, $d_2$, $d_3$ into \mythetag{3.7}, 
we get the identity $w_2=w_2$. However, there are the other three functions 
$d_1(b,c,w_1)$, $d_2(b,c,w_1)$, $d_3(b,c,w_1)$ in \mythetag{1.27}. They represent 
the same three roots $d_1$, $d_2$, $d_3$ of the same cubic equation \mythetag{1.19}. 
Substituting them into \mythetag{3.7}, we get the same quantity $w_2$. But now $w_2$ 
turns out to be expressed through $w_1$, i\.\,e\. we get the formula 
$$
\hskip -2em
w_2=\frac{3\,d_1(b,c,w_1)-3\,d_2(b,c,w_1)}{2\,d_3(b,c,w_1)-d_1(b,c,w_1)-d_2(b,c,w_1)},
\mytag{3.8}
$$
which is not an identity. The formula \mythetag{3.8} is the second conversion formula. 
It expresses $w_2$ through $b$, $c$, and $w_1$, i\.\,e\. \mythetag{3.8} yields a 
function $w_2=w_1(b,c,w_1)$. If we substitute this function into the argument
$w_2$ of the functions $d_1(b,c,w_2)$, $x_2(b,c,w_2)$, $x_3(b,c,w_2)$ from 
\mythetag{1.28}, then, according to the theorem~\mythetheorem{3.1}, we get back
three roots $d_1(b,c,w_1)$, $d_2(b,c,w_1)$, $d_3(b,c,w_1)$ used in \mythetag{3.8}.
This means that we have the following relationships based on the formula
\mythetag{3.8}:
$$
\align
&\hskip -2em
d_1(b,c,w_2(b,c,w_1))=d_1(b,c,w_1),\\
&\hskip -2em
d_2(b,c,w_2(b,c,w_1))=d_2(b,c,w_1),
\mytag{3.9}\\
&\hskip -2em
d_3(b,c,w_2(b,c,w_1))=d_3(b,c,w_1).
\endalign	
$$\par
     The formulas \mythetag{3.5} and \mythetag{3.8} provide two transformations
$w_1=w_1(b,c,w_2)$ and $w_2=w_2(b,c,w_1)$. Our next step is to prove that these
transformations are inverse to each other. For this purpose let's recall that
the functions \mythetag{1.27} obey the relationships \mythetag{1.5}, 
\mythetag{1.6}, \mythetag{1.7} whose right hand sides are given by the formulas
\mythetag{1.9} through \mythetag{1.17}. \pagebreak The same is true for the 
functions \mythetag{1.28}. In particular, we have the following three relationships 
for the functions \mythetag{1.28}:
$$
\aligned
&x_1(b,c,w_2)\,x_2(b,c,w_2)\,d_3(b,c,w_2)+x_2(b,c,w_2)\,x_3(b,c,w_2)\,\times\\
&\quad\times\,d_1(b,c,w_2)+x_3(b,c,w_2)\,x_1(b,c,w_2)\,d_2(b,c,w_2)=E_{21}(b,c),\\
\vspace{1ex}
&x_1(b,c,w_2)\,d_2(b,c,w_2)+d_1(b,c,w_2)\,x_2(b,c,w_2)+x_2(b,c,w_2)\,\times\\
&\quad\times\,d_3(b,c,w_2)+d_2(b,c,w_2)\,x_3(b,c,w_2)+x_3(b,c,w_2)\,\times\\
&\quad\times d_1(b,c,w_2)+d_3(b,c,w_2)\,x_1(b,c,w_2)=E_{11}(b,c),\\
\vspace{1ex}
&d_1(b,c,w_2)+d_2(b,c,w_2)+d_3(b,c,w_2)=E_{01}(b,c).
\endaligned
\mytag{3.10}
$$
Let's apply the formulas \mythetag{3.6} to \mythetag{3.10}. This yields
$$
\aligned
&x_1(b,c,w_1)\,x_2(b,c,w_1)\,d_3+x_2(b,c,w_1)\,x_3(b,c,w_1)\,d_1\,+\\
&\kern 10em +\,x_3(b,c,w_1)\,x_1(b,c,w_1)\,d_2=E_{21}(b,c),\\
\vspace{1ex}
&x_1(b,c,w_1)\,d_2+d_1\,x_2(b,c,w_1)+x_2(b,c,w_1)\,d_3\,+\\
&\ +\,d_2\,x_3(b,c,w_1)+x_3(b,c,w_1)\,d_1+d_3\,x_1(b,c,w_1)=E_{11}(b,c),\\
\vspace{1ex}
&d_1+d_2+d_3=E_{01}(b,c),
\endaligned
\mytag{3.11}
$$
where $w_1=w_1(b,c,w_2)$ and $d_i=d_i(b,c,w_2)$. The equalities \mythetag{3.11} 
are linear with respect to $d_1$, $d_2$, $d_3$. They constitute that very system 
of linear equations which was used in deriving the functions $d_i=d_i(b,c,w_1)$
(see \thetag{3.5} in \mycite{75}). This yields
$$
\align
&\hskip -2em
d_1(b,c,w_1(b,c,w_2))=d_1(b,c,w_2),\\
&\hskip -2em
d_2(b,c,w_1(b,c,w_2))=d_2(b,c,w_2),
\mytag{3.12}\\
&\hskip -2em
d_3(b,c,w_1(b,c,w_2))=d_3(b,c,w_2).
\endalign	
$$\par
    Apart from \mythetag{3.10} one can extract other three equations from
\mythetag{1.5}, \mythetag{1.6}, \mythetag{1.7} and write them as equalities 
for the functions \mythetag{1.27}, i\.\,g\. we can write
$$
\aligned
&x_1(b,c,w_1)+x_2(b,c,w_1)+x_3(b,c,w_1)=E_{10}(b,c),\\
\vspace{1ex}
&x_1(b,c,w_1)\,d_2(b,c,w_1)+d_1(b,c,w_1)\,x_2(b,c,w_1)+x_2(b,c,w_1)\,\times\\
&\quad\times\,d_3(b,c,w_1)+d_2(b,c,w_1)\,x_3(b,c,w_1)+x_3(b,c,w_1)\,\times\\
&\quad\times d_1(b,c,w_1)+d_3(b,c,w_1)\,x_1(b,c,w_1)=E_{11}(b,c),\\
\vspace{1ex}
&x_1(b,c,w_1)\,d_2(b,c,w_1)\,d_3(b,c,w_1)+x_2(b,c,w_1)\,d_3(b,c,w_1)\,\times\\
&\quad\times\,d_1(b,c,w_1)
+x_3(b,c,w_1)\,d_1(b,c,w_1)\,d_2(b,c,w_1)=E_{12}(b,c).
\endaligned
\mytag{3.13}
$$
Like in the case of \mythetag{3.10}, applying \mythetag{3.9} to \mythetag{3.13},
we get 
$$
\aligned
&x_1+x_2+x_3=E_{10}(b,c),\\
\vspace{1ex}
&x_1\,d_2(b,c,w_2)+d_1(b,c,w_2)\,x_2+x_2\,d_3(b,c,w_2)\,+\\
&\ +\,d_2(b,c,w_2)\,x_3+x_3\,d_1(b,c,w_2)+d_3(b,c,w_2)\,x_1=E_{11}(b,c),\\
\vspace{1ex}
&x_1\,d_2(b,c,w_2)\,d_3(b,c,w_2)+x_2\,d_3(b,c,w_2)\,d_1(b,c,w_2)\,+\\
&\kern 10em +\,x_3\,d_1(b,c,w_2)\,d_2(b,c,w_2)=E_{12}(b,c),
\endaligned
\mytag{3.14}
$$
where $w_2=w_2(b,c,w_1)$ and $x_i=x_i(b,c,w_1)$. The equalities \mythetag{3.14} 
are linear with respect to $x_1$, $x_2$, $x_3$. They constitute that very system 
of linear equations which was used in deriving the functions $x_i=x_i(b,c,w_2)$
(see \thetag{4.5} in \mycite{75}). This yields
$$
\align
&\hskip -2em
x_1(b,c,w_2(b,c,w_1))=x_1(b,c,w_1),\\
&\hskip -2em
x_2(b,c,w_2(b,c,w_1))=x_2(b,c,w_1),
\mytag{3.15}\\
&\hskip -2em
x_3(b,c,w_2(b,c,w_1))=x_3(b,c,w_1).
\endalign	
$$
The relationships \mythetag{3.15} are similar to \mythetag{3.6} and the 
relationships \mythetag{3.12} are similar to \mythetag{3.9}. But these four
groups of relationships do not coincide with each other.\par
     Now let's consider the composite function $w_2(b,c,w_1(b,c,w_2))$. Then, 
applying the formula \mythetag{3.8} to this function, we derive 
$$
\gathered
w_2(b,c,w_1(b,c,w_2))=\bigl(3\,d_1(b,c,w_1(b,c,w_2))
-3\,d_2(b,c,w_1(b,c,w_2))\bigr)\,\times\\
\times\,\bigl(2\,d_3(b,c,w_1(b,c,w_2))-d_1(b,c,w_1(b,c,w_2))
-d_2(b,c,w_1(b,c,w_2))\bigr)^{-1}.
\endgathered
$$
If we take into account \mythetag{3.12}, then the above formula can be 
written as
$$
\gathered
w_2(b,c,w_1(b,c,w_2))=\frac{3\,d_1(b,c,w_2)
-3\,d_2(b,c,w_2)}{2\,d_3(b,c,w_2)-d_1(b,c,w_2)-d_2(b,c,w_2)}.
\quad
\endgathered
\mytag{3.16}
$$
The right hand side of the formula \mythetag{3.16} can be produced by 
substituting the functions $d_1(b,c,w_2)$, $d_2(b,c,w_2)$, $d_3(b,c,w_2)$
from \mythetag{1.28} for $d_1$, $d_2$, $d_3$ into the formula
\mythetag{3.7}. The formula \mythetag{3.7} is a version of the formula
\mythetag{3.1} for $w=w_2$, while $d_1(b,c,w_2)$, $d_2(b,c,w_2)$, $d_3(b,c,w_2)$ 
are the roots of the cubic equation \mythetag{1.19} produced by means of the
formulas \mythetag{2.3}, \mythetag{2.4}, \mythetag{2.5} with $w=w_2$ applied
to the cubic equation \mythetag{1.19}. Therefore, the theorem~\mythetheorem{3.1}
in this case means that the right hand side of \mythetag{3.16} is equal to
$w_2$. Thus, we have derived the formula 
$$
\hskip -2em
w_2(b,c,w_1(b,c,w_2))=w_2.
\mytag{3.17}
$$\par
     The formula $w_1(b,c,w_2(b,c,w_1))=w_1$ is derived similarly. For this 
purpose we consider the function $w_1(b,c,w_2(b,c,w_1))$ and apply the formula 
\mythetag{3.5} to it:
$$
\gathered
w_1(b,c,w_2(b,c,w_1))=\bigl(3\,x_1(b,c,w_2(b,c,w_1))
-3\,x_2(b,c,w_2(b,c,w_1))\bigr)\,\times\\
\times\,\bigl(2\,x_3(b,c,w_2(b,c,w_1))-x_1(b,c,w_2(b,c,w_1))
-x_2(b,c,w_2(b,c,w_1))\bigr)^{-1}.
\endgathered
$$
Using the relationships \mythetag{3.15}, the above formula is transformed to 
$$
\gathered
w_1(b,c,w_2(b,c,w_1))=\frac{3\,x_1(b,c,w_1)
-3\,x_2(b,c,w_1)}{2\,x_3(b,c,w_1)-x_1(b,c,w_1)-x_2(b,c,w_1)}.
\quad
\endgathered
\mytag{3.18}
$$
The right hand side of the formula \mythetag{3.18} can be produced by 
substituting the functions $x_1(b,c,w_1)$, $x_2(b,c,w_1)$, $x_3(b,c,w_1)$
from \mythetag{1.27} for $x_1$, $x_2$, $x_3$ into the formula
\mythetag{3.1}, where $w=w_1$, while $x_1(b,c,w_2)$, $x_2(b,c,w_2)$, $x_3(b,c,w_2)$ 
are the roots of the cubic equation \mythetag{1.18} produced by means of the
formulas \mythetag{2.3}, \mythetag{2.4}, \mythetag{2.5} with $w=w_1$ applied
to the cubic equation \mythetag{1.18}. Therefore, the theorem~\mythetheorem{3.1}
in this case means that the right hand side of \mythetag{3.18} is equal to
$w_1$. Thus, we have derived the required formula for the composite function
$w_1(b,c,w_2(b,c,w_1))$:
$$
\hskip -2em
w_1(b,c,w_2(b,c,w_1))=w_1.
\mytag{3.19}
$$\par
     The functions $w_1=w_1(b,c,w_2)$ and $w_2=w_2(b,c,w_1)$ are given by
the formulas \mythetag{3.5} and \mythetag{3.8}. Using the explicit formulas 
for the functions \mythetag{1.27} and \mythetag{1.28} given in the
ancillary files {\bf Solutions\_\kern 1.5pt 1.txt} and {\bf Solutions\_\kern 1.5pt 
2.txt}, the formulas \mythetag{3.5} and \mythetag{3.8} are converted to explicit 
formulas for the functions $w_1=w_1(b,c,w_2)$ and $w_2=w_2(b,c,w_1)$. These 
explicit formulas are placed in the ancillary file {\bf Conversion\_\kern 1.5pt 
formulas.txt} attached to this arXiv submission.\par
      Theoretically, we could prove the relationships \mythetag{3.17} and
\mythetag{3.19} by direct calculations. However, the explicit formulas for
the functions $w_1(b,c,w_2)$ and $w_2(b,c,w_1)$ are so huge that they cannot 
be handled on my personal computer.\par
\head
4. Conclusions. 
\endhead
     The formulas \mythetag{3.17} and \mythetag{3.19} mean that the transformations
given by the conversion functions \mythetag{3.5} and \mythetag{3.8} are inverse to 
each other. Then the formulas \mythetag{3.6}, \mythetag{3.9}, \mythetag{3.12}, and 
\mythetag{3.15} mean that \mythetag{1.27} and \mythetag{1.28} are not two different
solutions of the cuboid factor equations \mythetag{1.3} and \mythetag{1.4}, but
two presentations of a single solution. This fact is the main result of the present
paper.\par

\enddocument
\end